\documentclass[11pt]{article}
\usepackage[T1]{fontenc}
\usepackage[utf8]{inputenc}
\usepackage{lmodern,amsmath,amssymb,amsthm,mathtools,booktabs,microtype}
\usepackage[a4paper,margin=27mm]{geometry}
\usepackage[colorlinks=true,linkcolor=blue,citecolor=blue,urlcolor=blue]{hyperref}
\hypersetup{pdftitle={Explicit exponent bounds for sums of consecutive powers},
  pdfauthor={Ákos Pintér},pdfdisplaydoctitle=true}
\usepackage{enumitem}
\usepackage{placeins}
\setlist{nosep}
\newtheorem{theorem}{Theorem}[section]
\newtheorem{lemma}[theorem]{Lemma}
\newtheorem{proposition}[theorem]{Proposition}
\newtheorem{corollary}[theorem]{Corollary}
\theoremstyle{remark}\newtheorem{remark}[theorem]{Remark}
\newcommand{\Z}{\mathbb Z}
\title{Explicit exponent bounds\\for sums of consecutive powers}
\author{\'{A}kos Pint\'{e}r\\[3pt]
\small Institute of Mathematics, University of Debrecen\\
\small 4032 Debrecen, Egyetem t\'{e}r 1, Hungary\\
\small\href{mailto:apinter@science.unideb.hu}{\texttt{apinter@science.unideb.hu}}}
\date{}
\begin{document}
\maketitle
\begin{abstract}
Let $x^k+(x+1)^k=y^n$, where $x,y,k,n$ are positive integers with
$k\ge2$ and $n\ge3$. We prove that
$n<52k^2\log(2k)\log^2(50k)$ for odd $k$ and
$n<216k\log^2(2500k)$ for even $k$. For odd $k$, an exact valuation
formula reduces the problem to a linear form in two real logarithms.
For even $k$, we use factorization in $\mathbb Z[i]$ and a bound of Yamada.
When $k=5$, we prove that the original exponent $n$ is one of the
34 primes between $23$ and $181$, inclusive. Every solution satisfies
$2x+1=u^n$ with $u>10^{1000}$. A gap between successive rational
approximations to $2^{-4/n}$, together with the Thue--Siegel principle,
gives at most four solutions for each of $n=23,29$ and at most three
for each larger prime in this range. Thus there are at most 104 positive
integer triples $(x,y,n)$ in the quintic case. All bounds are unconditional.
\end{abstract}



\section{Introduction}\label{sec:introduction}

For a fixed integer $k\ge2$, consider the equation
\begin{equation}\label{main}
 x^k+(x+1)^k=y^n,\qquad x\ge1,\quad y>1,\quad n\ge3,
\end{equation}
in integers $x,y,n$. Schinzel and Tijdeman \cite{ST} proved that the
exponent in $P(x)=y^n$ is effectively bounded when $P\in\mathbb Q[X]$
has at least two distinct zeros and $|y|>1$. Here
$P_k(X)=X^k+(X+1)^k$ has the $k$ distinct zeros
$\zeta/(1-\zeta)$ as $\zeta$ ranges over the roots of $Z^k=-1$. Hence the effective finiteness of all triples
$(x,y,n)$ for fixed $k$ also follows from their first corollary.
Our first results give explicit bounds in terms of $k$,
of order $k^2(\log k)^3$ for odd $k$ and $k(\log k)^2$ for even $k$.

Koutsianas and Tzanakis \cite[Theorems~1.1 and~1.3]{KT} exclude positive
solutions for $k\equiv2\pmod4$ whenever $6\le k\le100$. Their
result also covers every $k\equiv2\pmod4$ with $k\ge6$ whose odd prime
divisors are all congruent to $3$ modulo $4$. The case $4\mid k$
follows from results on the generalized Fermat equation
\cite{BEN,KT}. At $k=2$ there is a positive solution,
$119^2+120^2=13^4$. For comparison, at $k=50$ the estimate in
Theorem~\ref{even-bound} gives $n<994\,117$. The intermediate bound in
\cite[Table~2]{KT}, used there in the proof of nonexistence, is
$7\,627\,398$.

Our second result concerns $k=5$. Since $x$ and $x+1$ are coprime,
the known exclusions for the signatures $(5,5,p)$ imply that every
prime divisor of $n$ is at least $23$ \cite{DS,Poonen,Stoll}.
We apply the first logarithmic estimate to $n$ itself and obtain
$n<479<23^2$. It follows that $n$ is prime. The descent writes $2x+1=u^n$ and produces a rational approximation
$v/u^4$ to $2^{-4/n}$.
Continued fractions exclude $u\le10^{1000}$, and Laurent's determinant
estimate leaves precisely 34 candidate exponents. For each such prime,
two solutions give two particularly close approximations to the same
algebraic number. The gap between their denominators and the
Thue--Siegel principle together bound their number. The argument leaves the existence of a positive solution open for these 34 exponents.

Section~\ref{sec:theorems} states the results. Section~\ref{sec:general}
proves the general bounds and the extension to a fixed odd step.
Sections~\ref{sec:quintic} and~\ref{sec:counting} give the exponent
reduction and the solution count for $k=5$. The verification programs
are described in Section~\ref{sec:computations}. All logarithms are natural.

\section{Main results}\label{sec:theorems}

\subsection*{Bounds for general $k$}
For a prime $q$ and a positive integer $a$, let $v_q(a)$ be the exponent
of $q$ in $a$. For positive integers $m,A$, define
\[
 D_m(A)=\prod_{q\mid\gcd(m,A)}q^{v_q(m)}.
\]
Thus $D_m(A)$ contains the full $q$-primary part of $m$ for each prime
$q$ dividing $A$. In particular, $D_m(A)$ is a \emph{unitary divisor} of $m$:
\[
 D_m(A)\mid m,\qquad \gcd(D_m(A),m/D_m(A))=1.
\]
The empty product is $1$. We write $\omega(m)$ for the number of
distinct prime divisors of $m$. For a positive rational $a/b$ in lowest
terms, its logarithmic height is $h(a/b)=\log\max\{a,b\}$.

For odd $k$ we give three forms of the exponent bound. The first retains
the divisor supplied by the descent. The last depends only on $k$.

\begin{theorem}\label{odd-bound}
For any solution of \eqref{main} with odd $k\ge3$, put $d=D_k(2x+1)$ and
\[
 H=(k-1)\log2+k\log d,\qquad
 M(n)=\max\left\{\log\left(\frac nH+
 \frac1{(k-1)\log3}\right)+0.38,\ 10\right\}.
\]
Then
\begin{equation}\label{general-bound}
 n<12.6(k-1)H M(n)^2+
 \frac{\log\!\left(k(k-1)d^2/2\right)}{2\log3}.
\end{equation}
For each fixed $k$, this is an effective upper bound for $n$, with at
most $2^{\omega(k)}$ choices of $d$.
\end{theorem}

Replacing the unitary divisor by a quantity determined by $k$ gives
the following uniform version.

\begin{corollary}[Eliminating the unitary divisor]\label{k-only-bound}
Let $k\ge3$ be odd and put
\[
 H_k=(k-1)\log2+k\log k,\qquad
 M_k(n)=\max\left\{\log\left(\frac n{H_k}
       +\frac1{(k-1)\log3}\right)+0.38,\ 10\right\}.
\]
Every solution of \eqref{main} satisfies
\begin{equation}\label{bound-only-k}
 n<12.6(k-1)H_k M_k(n)^2+
 \frac{\log\!\left(k^3(k-1)/2\right)}{2\log3}.
\end{equation}
This gives an effective upper bound for $n$ in terms of $k$ alone,
without enumerating the unitary divisors of $k$.
\end{corollary}

The occurrence of $n$ in this estimate can also be eliminated explicitly.

\begin{corollary}[A fully explicit bound in terms of $k$]\label{explicit-k-bound}
Let $k\ge3$ be odd, and define
\[
 H_k=(k-1)\log2+k\log k,\qquad
 L_k=\max\{5,\ \log(26(k-1))+0.38\}.
\]
Every solution of \eqref{main} satisfies
\begin{equation}\label{explicit-k-sharp}
 n<52(k-1)H_k L_k^2.
\end{equation}
In particular,
\begin{equation}\label{explicit-k-simple}
 n<52k^2\log(2k)\log^2(50k).
\end{equation}
In particular, the original exponent satisfies $n=O(k^2(\log k)^3)$
with an absolute implied constant.
\end{corollary}

The same logarithmic argument applies when the two bases differ by
a fixed positive odd integer.

\begin{corollary}[A fixed odd step]\label{odd-step}
Let $h$ be a fixed positive odd integer, $x\ge1$, and $\gcd(x,h)=1$. For a solution of
\[
 x^k+(x+h)^k=y^n,\qquad k\ge3\text{ odd},\quad n\ge3,
\]
put $t=2x+h$, $d=D_k(t)$, and define $H,M(n)$ as in
Theorem~\ref{odd-bound} with this $d$. Then
\[
 n<12.6(k-1)H M(n)^2+
 \frac{\log\!\left(k(k-1)d^2h^2/2\right)}{2\log3}.
\]
\end{corollary}

For even $k$, Gaussian factorization gives a bound of smaller
asymptotic order through a logarithmic estimate for complex arguments.

\begin{theorem}\label{even-bound}
Let $k\ge2$ be even. Every solution of \eqref{main}, with arbitrary
integer $n\ge3$, satisfies
\begin{equation}\label{even-k-bound}
 n<216kL_k^2<216k\log^2(2500k),\qquad
 L_k=\max\left\{\frac{17}{2},\log(108k)+1\right\}.
\end{equation}
In particular $n=O(k(\log k)^2)$, with an absolute effective constant.
\end{theorem}

When $k\equiv2\pmod4$, the exact valuation descent also restricts
the coefficients of the auxiliary equations.

\begin{proposition}\label{even}
Let $k=2m\ge2$ with $m$ odd, and suppose \eqref{main} holds. Put
\[
 A=x^2+(x+1)^2,\qquad
 B=\frac{x^{2m}+(x+1)^{2m}}A,\qquad d=D_m(A).
\]
There exist positive odd integers $u,v$ satisfying
\begin{equation}\label{even-split}
 dA=u^n,\qquad B=dv^n,\qquad y=uv,\qquad\gcd(u,v)=1.
\end{equation}
The coefficient $d$ is a unitary divisor of $m$, and every prime dividing
$d$ is congruent to $1\pmod4$.
\end{proposition}

If every possible coefficient is itself a $p$th power, this descent
reduces to a known equation and yields the following obstruction.

\begin{corollary}\label{valuation-obstruction}
Let $p\ge5$ be prime and $k=2m$ with $m$ odd. If
$p\mid v_q(m)$ for every prime $q\equiv1\pmod4$ dividing $m$, then
$x^k+(x+1)^k=y^p$ has no solution with $x\ge1$.
\end{corollary}

\subsection*{The quintic equation}

\begin{theorem}\label{quintic-main}
Let $x,y,n$ be positive integers with $x^5+(x+1)^5=y^n$ and $n\ge3$.
Then
\begin{equation}\label{E}
 n\in\mathcal E=\{p:p\text{ is prime},\ 23\le p\le181\}.
\end{equation}

There are positive odd integers $u,v$ such that
\begin{equation}\label{quintic-system}
 2x+1=u^n,\qquad y=uv,\qquad
 u^{4n}+10u^{2n}+5=16v^n,
\end{equation}
with $\gcd(u,v)=\gcd(uv,5n)=1$ and
\begin{equation}\label{huge}
u>10^{1000},\qquad x>\frac{10^{23000}-1}{2}.
\end{equation}
\end{theorem}

To see the converse, put $t=2x+1$. The identity
$16(x^5+(x+1)^5)=t(t^4+10t^2+5)$ shows that any positive odd
$u\ge3$ and positive $v$ satisfying the last equation of
\eqref{quintic-system} yield a solution with $x=(u^n-1)/2$ and $y=uv$.
For fixed $n$, each such $u$ determines $v$ uniquely.

\begin{proposition}\label{gap}
Fix $n\in\mathcal E$. If two distinct quintic solutions have auxiliary
variables $u_1<u_2$, then
\[
 u_2>u_1^{(n-2)/2}.
\]
Equivalently their original variables satisfy
$2x_2+1>(2x_1+1)^{(n-2)/2}$.
In particular, a second solution for the same exponent must have
$u_2>10^{500(n-2)}\ge10^{10500}$.
\end{proposition}

\begin{theorem}\label{quintic-count}
For $n\in\mathcal E$, let
\[
 N_n=\#\{(x,y)\in\Z_{>0}^2:x^5+(x+1)^5=y^n\}.
\]
Then
\[
 N_n\le
 \begin{cases}
  4,& n=23\text{ or }29,\\
  3,& 31\le n\le181.
 \end{cases}
\]
Consequently, there are at most $104$ triples $(x,y,n)$ of positive
integers satisfying $x^5+(x+1)^5=y^n$ and $n\ge3$.
\end{theorem}

\begin{corollary}\label{quintic-distinct-x}
For each positive integer $x$, there is at most one pair $(y,n)$
with $y>1$, $n\ge3$ and $x^5+(x+1)^5=y^n$. In particular, at most
$104$ positive integers $x$ admit such a representation.
\end{corollary}

\section{Descent and general exponent bounds}\label{sec:general}
\subsection{The valuation descent}

\begin{lemma}\label{descent}
Let $a,b$ be coprime positive integers of opposite parity, and let $m$
be a positive odd integer. Put
\[
 A=a+b,\qquad B=\frac{a^m+b^m}{a+b},\qquad d=D_m(A).
\]
Then $B/d$ is a positive integer and
\begin{equation}\label{coprime-split}
 \gcd(dA,B/d)=1.
\end{equation}
Consequently, for any integer $n\ge2$, the equation $a^m+b^m=y^n$
is equivalent to
\begin{equation}\label{split}
 dA=u^n,\qquad B=dv^n,\qquad y=uv
\end{equation}
for coprime positive odd integers $u,v$.
\end{lemma}

\begin{proof}[Proof of Lemma~\ref{descent}]
Since $m$ is odd, $A$ divides $a^m+b^m$. The two bases have opposite
parity, so $A$ is odd. Any prime $q\mid A$ divides neither $a$ nor $b$. The lifting-the-exponent
identity gives
\[
 v_q(a^m+b^m)=v_q(A)+v_q(m).
\]
Hence $v_q(B)=v_q(m)$ for every $q\mid A$. All such prime factors are
removed completely from $B$ upon division by $d$, proving
\eqref{coprime-split}. Since
\[
 a^m+b^m=(dA)(B/d),
\]
the product is an $n$th power precisely when its two coprime positive
factors are $n$th powers. Their oddness is immediate.
\end{proof}

\begin{remark}
The lemma has no hypothesis comparing $n$ with $v_q(m)$. A prime dividing
$m$ is either absent from $d$ or occurs in $d$ to its full exponent in $m$.
Thus there are at most $2^{\omega(m)}$ possible coefficients, where
$\omega(m)$ is the number of distinct prime divisors of $m$.
\end{remark}

\subsection{Odd powers}\label{sec:odd}

Let $k\ge3$ be odd and set $t=2x+1$. Define
\[
 G_k(T)=\sum_{j=0}^{(k-1)/2}\binom{k}{2j}T^{k-1-2j}.
\]
Then
\begin{equation}\label{identity}
 x^k+(x+1)^k=\frac{tG_k(t)}{2^{k-1}}.
\end{equation}
Applying Lemma~\ref{descent} with $(a,b,m)=(x,x+1,k)$ gives
\begin{equation}\label{odd-split}
 dt=u^n,\qquad G_k(t)=2^{k-1}dv^n,\qquad
 d=D_k(t),\quad \gcd(u,v)=1.
\end{equation}
Here $u,v$ are odd and $u\ge3$.

\begin{proof}[Proof of Theorem~\ref{odd-bound}]
Set $R=G_k(t)/t^{k-1}$. Since $t\ge3$,
\[
 1<R<2^{k-1}.
\]
Moreover, for each $j\ge0$ in the sum defining $G_k$,
\[
 \binom{k}{2j}\le \frac{(k(k-1)/2)^j}{j!}.
\]
Indeed, pair the factors of $k(k-1)\cdots(k-2j+1)$ and use
$(2j)!\ge2^j j!$. Consequently
\begin{equation}\label{log-R}
 0<\log R<\frac{k(k-1)}{2t^2}.
\end{equation}
With
\[
 \alpha=\frac{u^{k-1}}v,\qquad \beta=2^{k-1}d^k,
\]
equations \eqref{odd-split} give $\alpha^n=\beta/R>1$ and
\begin{equation}\label{Lambda-general}
 \Lambda:=\log\beta-n\log\alpha=\log R>0.
\end{equation}
The number $\alpha$ is a ratio of odd integers, whereas $v_2(\beta)=k-1>0$.
Since $\alpha>1$, these two rational numbers are multiplicatively
independent. Also, $v<u^{k-1}$ and $\gcd(u,v)=1$, so that
\[
 h(\alpha)=(k-1)\log u,\qquad h(\beta)=H.
\]
Laurent's Corollary~2 \cite{Laurent}, with $D=1$, $m=10$, $C_2=25.2$,
$b_1=n$, $b_2=1$, $\log A_1=(k-1)\log u$, and $\log A_2=H$, gives
\[
 \log\Lambda\ge-25.2(k-1)(\log u)H M(n)^2.
\]
On the other hand, $t=u^n/d$ and \eqref{log-R} imply
\[
 \log\Lambda<\log\!\left(k(k-1)d^2/2\right)-2n\log u.
\]
Combining these estimates and using $u\ge3$ proves \eqref{general-bound}.
Its right-hand side grows as $O_k((\log n)^2)$, and $d$ ranges over a
finite effectively determined set.
\end{proof}

\begin{proof}[Proof of Corollary~\ref{k-only-bound}]
Since $d=D_k(2x+1)$ is a unitary divisor of $k$, we have
$1\le d\le k$ and $H\le H_k$. The dependence of $M(n)$ on $H$
requires checking the whole product. Put $c=1/((k-1)\log3)$ and,
for $s>0$, define
\[
 F(s)=s\max\{\log(n/s+c)+0.38,10\}^2.
\]
On the constant branch $F(s)=100s$. On the other branch, writing
$L=\log(n/s+c)+0.38\ge10$, differentiation gives
\[
 F'(s)=L\left(L-\frac{2n}{n+cs}\right)>0.
\]
The two expressions agree at their common endpoint, so $F$ is
increasing on $(0,\infty)$. Hence
$H M(n)^2\le H_k M_k(n)^2$. Also
$k(k-1)d^2/2\le k^3(k-1)/2$.
Substitution in \eqref{general-bound} proves the claim.
\end{proof}

We use the following elementary estimate to eliminate the exponent from
both logarithmic bounds.

\begin{lemma}\label{invert-log}
Let $A>0$, $b\in\mathbb R$, $c>2$ and $W(t)=\max\{c,\log t+b\}$ for
$t>0$. Put $L=\max\{5,c/2,\log(2A)+b\}$. If $t<AW(t)^2$, then
$t<4AL^2$.
\end{lemma}
\begin{proof}
The function $t/W(t)^2$ is strictly increasing. Its derivative on the logarithmic
branch is $(W(t)-2)/W(t)^3>0$. Put $T=4AL^2$. Since $L\ge5$,
$\log2+2\log L<L$ (the difference $L-\log2-2\log L$ is positive at
$5$ and increasing thereafter). Thus $\log T+b<2L$, and $c\le2L$
gives $W(T)\le2L$. Now $t/W(t)^2<A\le T/W(T)^2$, so $t<T$.
\end{proof}

\begin{proof}[Proof of Corollary~\ref{explicit-k-bound}]
Put
\[
 \delta=\frac1{(k-1)\log3},\quad
 C_k=\frac{\log(k^3(k-1)/2)}{2\log3},\quad
 z=\frac n{H_k}+\delta,\quad
 W(z)=\max\{\log z+0.38,10\}.
\]
Since $H_k\ge k\log k$ and $k^3(k-1)/2<k^4$, we have
\[
 \frac{C_k}{H_k}+\delta
 <\frac1{\log3}\left(\frac2k+\frac1{k-1}\right)
 \le\frac7{6\log3}<2.
\]
Consequently, \eqref{bound-only-k} implies
\[
 z<12.6(k-1)W(z)^2+2<13(k-1)W(z)^2,
\]
where the last inequality uses $k-1\ge2$ and $W(z)\ge10$.

Lemma~\ref{invert-log} with $A=13(k-1)$, $b=0.38$ and $c=10$
gives $z<52(k-1)L_k^2$. Since $n<H_kz$, this proves
\eqref{explicit-k-sharp}.

Finally, $H_k<k\log(2k)$. The elementary inequalities
$\log150>5$ and $\log(50/26)>0.38$ show that
$L_k<\log(50k)$ for every $k\ge3$. Using also $k-1<k$
gives \eqref{explicit-k-simple}.
\end{proof}

\begin{proof}[Proof of Corollary~\ref{odd-step}]
The integers $x,x+h$ are coprime and have opposite parity, so
Lemma~\ref{descent} again gives $dt=u^n$. Writing
$x=(t-h)/2$ and $x+h=(t+h)/2$ replaces $R$ by
\[
 R_h=\sum_{j=0}^{(k-1)/2}\binom{k}{2j}(h/t)^{2j}.
\]
As $0<h/t<1$, we have $1<R_h<2^{k-1}$, and the coefficient estimate
in the proof of Theorem~\ref{odd-bound} gives
$\log R_h<k(k-1)h^2/(2t^2)$. The heights of $\alpha,\beta$ are as
before, so substitution of $t=u^n/d$ gives the asserted bound.
\end{proof}

\subsection{Even powers}\label{sec:even}

For the even-power equation we use the following bound for a linear
form involving the argument of an algebraic number of absolute value one.

\begin{lemma}[Yamada {\cite[Theorem~1.2]{Yamada}}]\label{yamada}
Let $\alpha$ be an algebraic number with $|\alpha|=1$ that is not a
root of unity, and let $b_1,b_2$ be nonzero integers. Using the principal
value of $\log\alpha$, put
\[
 \Lambda=b_2\log\alpha-\frac{b_1\pi i}{2}.
\]
Let $h(\alpha)$ denote the absolute logarithmic height of $\alpha$, and set
\[
 D=\frac{[\mathbb Q(\alpha):\mathbb Q]}2,\qquad
 a=9\pi+2D h(\alpha),\qquad
 b'=\frac{|b_1|}{a}+\frac{|b_2|}{9\pi},
\]
\[
 H=\max\{17,\ D,\ D(\log b'+2.96)+0.01\}.
\]
Then
\begin{equation}\label{yamada-bound}
 \log|\Lambda|>-(2.76701a+0.12945)H^2.
\end{equation}
\end{lemma}

\begin{proof}[Proof of Theorem~\ref{even-bound}]
Write $k=2r$. The integer $y$ is odd. Also $3\nmid y$, since the
sum of two consecutive even powers is congruent to $1$ or $2$ modulo
$3$. Thus $y\ge5$.

Put $\beta=(x+1)^r+ix^r$. Its norm is $y^n$, and $\beta$ and
$\overline\beta$ are coprime in $\Z[i]$. A common Gaussian prime
would divide both $2(x+1)^r$ and $2ix^r$, hence divide $2$, whereas
the norm is odd. Unique factorization therefore gives
\[
 \beta=\zeta\gamma^n,\qquad
 \zeta\in\{\pm1,\pm i\},\qquad
 \gamma=c+di\in\Z[i],\qquad c^2+d^2=y.
\]
Moreover $\gamma$ and $\overline\gamma$ are coprime, so
$\gcd(c,d)=1$ and $c,d$ have opposite parity. Both coordinates
are nonzero, because both coordinates of $\beta$ are nonzero.
Consequently
\[
 \alpha=\frac{\gamma}{\overline\gamma}\in\mathbb Q(i)
\]
has absolute value one and is not a root of unity. Indeed, the four
roots of unity in $\mathbb Q(i)$ would require $cd=0$ or $c=\pm d$.
Conjugation sends $\alpha$ to $\alpha^{-1}$. The polynomial
\[
 yZ^2-2(c^2-d^2)Z+y,
\]
is primitive. An odd prime dividing both $y=c^2+d^2$ and $c^2-d^2$
would divide both $c$ and $d$, and $y$ is odd. Since $\alpha$ is
nonreal, this is its minimal polynomial over $\mathbb Z$ up to sign.
Both roots have modulus one, so $h(\alpha)=\tfrac12\log y$.

Let $\phi=\arg\alpha\in(-\pi,\pi]$ and
\[
 \theta=2\arctan\!\left(\left(\frac{x}{x+1}\right)^r\right)
 \in(0,\pi/2).
\]
Since $\beta/\overline\beta=\zeta^2\alpha^n=e^{i\theta}$,
we have $n\phi=\theta+j\pi$ for some $j\in\Z$. Set $b=2j+1$.
Using the principal logarithm, we obtain the nonzero form
\begin{equation}\label{even-form}
 \Lambda=n\log\alpha-\frac{b\pi i}{2}
       =i(\theta-\pi/2),\qquad 0<|b|<2n+1.
\end{equation}
For $s=r\log(1+1/x)>0$, elementary differentiation gives
\[
 0<|\Lambda|=\frac\pi2-2\arctan(e^{-s})
       =\int_0^s\frac{dt}{\cosh t}<s<\frac{k}{2x}.
\]
In particular,
\begin{equation}\label{even-upper}
 \log|\Lambda|<\log k-\log x.
\end{equation}

Apply Lemma~\ref{yamada} to \eqref{even-form}, with $b_1=b$ and $b_2=n$.
Its parameters are $D=[\mathbb Q(\alpha):\mathbb Q]/2=1$ and
\[
 a=9\pi+\log y,\qquad
 b'=\frac{|b|}{a}+\frac{n}{9\pi},\qquad
 H=\max\{17,\log b'+2.97\}.
\]
Inequality \eqref{yamada-bound} gives
\[
 \log|\Lambda|>-(2.76701a+0.12945)H^2.
\]
Since $n\ge3$,
\[
 b'<\frac{3n+1}{9\pi}\le\frac{10n}{27\pi},\qquad
 H\le W(n):=\max\{17,\log n+1\}.
\]
The latter inequality uses $\log(10/(27\pi))+2.97<1$.
Also $y\ge5$, $\pi<22/7$ and $\log5>8/5$ give
\[
 2.76701(9\pi+\log y)+0.12945<52\log y.
\]
We conclude that
\begin{equation}\label{even-lower}
 \log|\Lambda|>-52(\log y)W(n)^2.
\end{equation}

First suppose $x\le(2k)^{100}$. Then
\begin{align*}
 n\log5&\le n\log y<k\log(x+1)+\log2\\
 &\le100k\log(2k)+(k+1)\log2\le101k\log(2k),
\end{align*}
so $n<64k\log(2k)<216kL_k^2$.
Otherwise put $w=\log x>100\log(2k)$. We have
\[
 \log|\Lambda|<-0.99w,\qquad
 n\log y<kw+(k+1)\log2<1.01kw.
\]
Together with \eqref{even-lower}, these inequalities imply
\begin{equation}\label{even-implicit}
 n<52\frac{1.01}{0.99}kW(n)^2<54kW(n)^2.
\end{equation}

Applying Lemma~\ref{invert-log} with $A=54k$, $b=1$ and $c=17$ gives
$n<216kL_k^2$. Finally, $\log5000>17/2$ and
$\log(108k)+1<\log(2500k)$ show that $L_k<\log(2500k)$
for $k\ge2$, proving \eqref{even-k-bound}.
\end{proof}

\begin{proof}[Proof of Proposition~\ref{even}]
Apply Lemma~\ref{descent} to $(a,b)=(x^2,(x+1)^2)$.
If $q\mid A$, then $q$ is odd and
$(x/(x+1))^2\equiv-1\pmod q$, so $q\equiv1\pmod4$.
\end{proof}

Writing $t=2x+1$ gives $d(t^2+1)=2u^n$. In the notation of
\cite[Lemma~4.2]{KT}, the corresponding representation can always be
chosen with $d_1=d$, $d_2=1$, with $d$ unitary in $k/2$.
For instance, the descent gives the following candidate coefficients.
\[
\begin{array}{c|ccccc}
 k&10&26&34&50&130\\ \hline
 d&1,5&1,13&1,17&1,25&1,5,13,65.
\end{array}
\]
In particular, $d=5$ is not needed when $k=50$.

\begin{proof}[Proof of Corollary~\ref{valuation-obstruction}]
Every possible $d$ in Proposition~\ref{even} is a $p$th power. From
$dA=u^p$ it follows that $A=w^p$ for a positive integer $w$.
Hence $(2x+1)^2+1=2w^p$, which has no solution with $2x+1>1$ by
\cite[Theorem~1.2]{KT}. The argument also applies when $m$ has no prime divisor congruent to $1$ modulo $4$.
\end{proof}

\subsection{Exponents divisible by three}
For completeness we recall the following consequence of the generalized
Fermat equation, recorded in \cite[Remark~8]{BCDY}.

\begin{proposition}\label{three}
If $3\mid k$ and $k\ge6$, equation \eqref{main} has no solution.
More generally,
\[
 a^{3r}+b^{3s}=c^n,\qquad r,s\ge2,\quad n\ge3,
\]
has no solution in coprime nonzero integers $a,b,c$.
\end{proposition}

\begin{proof}[Proof of Proposition~\ref{three}]
Kraus's necessary condition, with the small-exponent cases included as in
\cite[Proposition~7]{BCDY}, states that a primitive nonzero solution of
$A^3+B^3=C^n$, $n\ge3$, must satisfy $v_2(AB)=1$.
For $A=a^r$, $B=b^s$, coprimality gives either $v_2(AB)=0$ or
$v_2(AB)\ge2$, a contradiction. Take $a=x$, $b=x+1$ and $r=s=k/3$
for the first assertion.
\end{proof}

\section{The exponent in the quintic equation}\label{sec:quintic}
\subsection{Prime divisors of the exponent and the quintic descent}

Let $(x,y,n)$ be a solution with $k=5$ and $n\ge3$. For each prime
$p\mid n$, the triple $(x,x+1,y^{n/p})$ is a nontrivial primitive
solution of $A^5+B^5=C^p$. The cases $p=2,3$ are excluded by
\cite{DM,Poonen}, and $p=5$ by Fermat's last theorem. Dahmen and
Siksek \cite[Theorem~1]{DS} exclude $p=7,19$. Stoll's theorem on
$x^5+y^5=z^p$ \cite[Section~8, Theorem~8.8 in the preprint]{Stoll}
also excludes $p=11,13,17$ unconditionally. Its extension through
$p=53$ assumes GRH for $p\ge23$. We use only the unconditional cases.
It follows that
\begin{equation}\label{n-coprime}
 p\mid n,\quad p\text{ prime}\quad\Longrightarrow\quad p\ge23.
\end{equation}

Billerey, Chen, Dieulefait and Freitas
\cite[Theorems~3 and~4]{BCDF} give $p\nmid C$ and $5\nmid C$ for
each such triple. Since $C=y^{n/p}$, this gives
\begin{equation}\label{local-y}
 \gcd(y,5n)=1.
\end{equation}
If $5\mid2x+1$, then $5\mid x^5+(x+1)^5$, contrary to
\eqref{local-y}. Thus $d=D_5(2x+1)=1$. Applying \eqref{odd-split}
with the original exponent $n$ gives \eqref{quintic-system}, with
$u\ge3$ odd and $\gcd(u,v)=\gcd(uv,5n)=1$.

\subsection{An upper bound for the original exponent}

Put $L=\log u$ and $\alpha=u^4/v$. We use the logarithmic form
\begin{equation}\label{quintic-Lambda}
 \Lambda=4\log2-n\log\alpha
 =\log(1+10u^{-2n}+5u^{-4n}),\qquad
 0<\Lambda<10u^{-2n}.
\end{equation}

Indeed, $1+10z+5z^2<\exp(10z)$ for $z=u^{-2n}>0$. Also $v<u^4$,
since $16s^2-(s^2+10s+5)=5(s-1)(3s+1)>0$ for $s=u^{2n}$.
Thus $0<\log\alpha<4\log2/n$ and $h(\alpha)=4L$.
The rational number $\alpha>1$ is a quotient of odd integers, so
$\alpha$ and $2$ are multiplicatively independent.

Laurent's Corollary~2 \cite{Laurent} applies with $\alpha_1=\alpha$, $\alpha_2=2$,
$b_1=n$, $b_2=4$, $\log A_1=4L$, and $\log A_2=1$. This gives
\begin{equation}\label{5041-bound}
 n<50.4\max\{\log(n+1/\log3)+0.38,10\}^2+
 \frac{\log10}{2\log3}.
\end{equation}
At $n=5042$ the right-hand side is less than $5041.048$. Dividing it
by $n$ gives a decreasing function for $n\ge5042$. Indeed, put
$\delta=1/\log3$. On the logarithmic branch,
$M=\log(n+\delta)+0.38\ge10$, and the derivative of $M^2/n$
has the sign of $2n/(n+\delta)-M<0$. On the constant branch the
claim is immediate. Thus
\begin{equation}\label{initial-n}
 n\le5041.
\end{equation}

To improve \eqref{initial-n} we specialize Laurent's Theorem~2
\cite{Laurent}. We keep the dependence on a lower bound for $u$
because the same estimate will be used again after the sieve and
the continued-fraction calculation. Suppose that
\[
 0<P\le n\le B,\qquad u\ge U>1,\qquad L_0=\log U.
\]
Choose $\rho>1$ and $1/3\le\mu\le1$, and define
\begin{align*}
 \sigma&=(1+2\mu-\mu^2)/2,& \lambda&=\sigma\log\rho,\\
 K&=8+\frac{4(\rho-1)\log2}{PL_0},& a_{10}&=KL_0,\\
 a_2&=(\rho+1)\log2,& H_0&=h/\lambda+1/\sigma,\\
 h&=\max\left\{\log\left(\frac B{a_2}+\frac4{a_{10}}\right)
       +\log\lambda+1.81,\lambda,\tfrac12\log2\right\},&&\\
 \omega&=2\left(1+\sqrt{1+\frac1{4H_0^2}}\right),&
 \theta&=\sqrt{1+\frac1{4H_0^2}}+\frac1{2H_0}.
\end{align*}
Let
\begin{align*}
 C_0={}&\frac{\mu}{\lambda^3\sigma}
 \left(\frac\omega6+\frac12
 \sqrt{\frac{\omega^2}9+
 \frac{8\lambda\omega^{5/4}\theta^{1/4}}
 {3\sqrt{a_{10}a_2}\sqrt{H_0}}+
 \frac43\left(\frac1{a_{10}}+\frac1{a_2}\right)
 \frac{\lambda\omega}{H_0}}\right)^2,\\
 C'_0={}&\sqrt{\frac{C_0\sigma\omega\theta}{\lambda^3\mu}},
 \qquad T=h+\lambda/\sigma,\qquad D_0=C'_0T^2Ka_2.
\end{align*}
In each application below the conditions
\begin{equation}\label{numerical-conditions}
 a_{10},a_2\ge1,\qquad a_{10}a_2\ge\lambda^2,
 \qquad \log(D_0L_0)>1
\end{equation}
are checked by outward-rounded interval arithmetic. Define
\begin{equation}\label{Q}
 \mathcal Q(P,B,U;\rho,\mu)=\frac12\left(
 C_0T^2Ka_2+
 \frac{\sqrt{\omega\theta}\,T+\log(D_0L_0)+\log10}{L_0}
 \right).
\end{equation}

\begin{lemma}\label{uniform-laurent}
Under the preceding hypotheses and \eqref{numerical-conditions}, every
quintic solution with $P\le n\le B$, $u\ge U$ satisfies
$n<\mathcal Q(P,B,U;\rho,\mu)$.
\end{lemma}

\begin{proof}[Proof of Lemma~\ref{uniform-laurent}]
In Laurent's Theorem~2 take $a_1=K\log u$. The inequality $0<\log\alpha<4\log2/P$ from \eqref{quintic-Lambda} gives $a_1\ge(\rho-1)\log\alpha+2h(\alpha)$.
The choice of $h$ is valid uniformly because $n\le B$ and $a_1\ge a_{10}$.
The coefficient $C$ in that theorem decreases as $a_1$ increases, so
$C\le C_0$ and $C'\le C'_0$. Its lower bound therefore yields
\[
 \log\Lambda\ge-C_0T^2Ka_2L-\sqrt{\omega\theta}\,T-\log(D_0L).
\]
Combine this with $\log\Lambda<\log10-2nL$.
The function $\log(D_0L)/L$ decreases on $L\ge L_0$ by the last
condition in \eqref{numerical-conditions}. The remaining
positive terms divided by $L$ also decrease. Substituting $L_0$ proves
\eqref{Q}.
\end{proof}

The certified parameters used in the proof of Theorem~\ref{quintic-main}
are listed in Table~\ref{bounds}.
Each printed bound is rounded upwards. Within a block, an inequality
$n<c$ permits the next integer upper bound $B=\lfloor c\rfloor$.
The last line of each block contradicts the assumed lower bound $P$.

\begin{table}[htbp]
\centering
\caption{Certified upper bounds for $\mathcal Q$.}\label{bounds}
\begin{tabular}{rrrrrr}
\toprule
$P$ & $B$ & $U$ & $\rho$ & $\mu$ & Upper bound\\
\midrule
479 & 5041 & 3 & 16 & $1/3$ & 712.232\\
479 & 712  & 3 & 16 & $1/3$ & 511.529\\
479 & 511  & 3 & 16 & $1/3$ & 480.835\\
479 & 480  & 3 & 16 & $1/3$ & 475.157\\
\midrule
331 & 467 & $10^6+1$ & 18 & $12/25$ & 348.645\\
331 & 348 & $10^6+1$ & 18 & $12/25$ & 328.261\\
\midrule
277 & 317 & $10^{1000}$ & 19 & $27/50$ & 282.880\\
277 & 282 & $10^{1000}$ & 19 & $27/50$ & 275.823\\
\bottomrule
\end{tabular}
\end{table}

If $n\ge479$, the first four rows of Table~\ref{bounds}, starting
from $n\le5041$, give successive upper bounds $712$, $511$ and $480$.
The fourth row gives $n<475.157$, a contradiction. Hence $n<479$. Every composite
integer whose prime divisors are all at least $23$ is at least
$23^2=529$. By \eqref{n-coprime}, the original $n$ is therefore
prime, and the largest prime below $479$ is $467$.

\subsection{The congruence sieve}

For each of the 83 prime exponents
\[
 23\le n\le467,
\]
we exclude $3\le u\le10^6$, with $u$ odd and $\gcd(u,5n)=1$.
For a prime $q=2rn+1$, put
\[
 \mathcal P_{n,q}=\{b^n:b\in\mathbb F_q\},\qquad
 \mathcal R_{n,q}=\{a\in\mathbb F_q:
 16^{-1}(a^{4n}+10a^{2n}+5)\in\mathcal P_{n,q}\}.
\]
Every solution of \eqref{quintic-system} must have
$u\bmod q\in\mathcal R_{n,q}$. The program starts with all admissible
$u$ in the stated interval and applies these tests for increasing $r$,
retaining only survivors. For each of the 83 prime exponents the survivor set
becomes empty with $r\le500$. The certificate records every modulus and
the survivor count after every step. All primality tests use trial division,
and all residue calculations are exact. The largest multiplier actually
used is $r=51$.

Thus $u\ge10^6+1$. The second block of Table~\ref{bounds} now excludes
$n\ge331$. Since $n$ is prime, the largest remaining possibility is
$317$.

\subsection{The continued-fraction bound}

Set $\gamma_n=2^{-4/n}$. Equation \eqref{quintic-system}, concavity of
$z^{1/n}$ and $n\ge23$, $u\ge3$ give
\begin{align}
 0<\frac v{u^4}-\gamma_n
 &=\gamma_n\left((1+10u^{-2n}+5u^{-4n})^{1/n}-1\right)\notag\\
 &\le\frac{10+5u^{-2n}}{n\,u^{2n}}
 <\frac1{u^{2n}}<\frac1{2u^8}.\label{CF-ineq}
\end{align}
Since $\gcd(v,u^4)=1$, Legendre's continued-fraction criterion shows
that $v/u^4$ is a convergent to $\gamma_n$.

For all 58 remaining prime exponents
\[
 23\le n\le317,
\]
the program enumerates all convergents with denominator at most
$10^{4000}$. A denominator relevant to \eqref{CF-ineq} must be the
fourth power of an odd integer $u$, with $\gcd(u,5n)=1$.
Among denominators $u^4$ with $u\ge3$ satisfying these conditions,
the only candidate is $80/81$ at $n=223$. Its numerator is even,
whereas $v$ must be odd. In all, $450\,795$ partial quotients are
certified, including the first convergent beyond the cutoff for each
exponent. Therefore
\begin{equation}\label{u-cf}
 u>10^{1000}.
\end{equation}

With $S=10^{8050}$ the program first constructs an integer $A$ and proves
by exact integer arithmetic that
\begin{equation}\label{alpha-cert}
 16A^n<S^n<16(A+1)^n.
\end{equation}
Thus $A/S<\gamma_n<(A+1)/S$. A partial quotient is accepted only when
both rational endpoints have the same integer part. After subtracting
that part and reciprocating, the endpoints are exchanged. This process
continues until a convergent denominator exceeds $10^{4000}$.
The certificate contains $A$, all accepted partial quotients, and the
first denominator beyond the cutoff. A decimal approximation is used
only to suggest $A$. Inequality \eqref{alpha-cert} certifies the interval.

The third block of Table~\ref{bounds}, using \eqref{u-cf}, excludes
$n\ge277$. Primality then gives $n\le271$.

\subsection{The final determinant estimate}\label{sec:direct}
We now assume $191\le n\le B\le271$ and choose the integer
parameters in Laurent's Theorem~1 as functions of $\log u$. This
makes the resulting bound uniform over $u>10^{1000}$.

Put $w=\log u$ and $w_0=1000\log10$, so $w>w_0$. Let $\kappa>0$ and let $\ell,r\ge2$ be integers, with
$r<191$, and take $\rho>1$ and $\mu=3/5$. Define
\begin{align*}
 \kappa_+&=\kappa+1/w_0,& K_0&=\kappa w_0,&
 s_0&=\ell/r+\ell/K_0,\\
 g_0&=\frac14-\frac{\ell}{12rs_0},&
 a&=8+\frac{4(\rho-1)\log2}{191w_0},& a_2&=(\rho+1)\log2,\\
 \sigma&=(1+2\mu-\mu^2)/2,&
 b_0&=\log\!\left(\frac{4r/K_0+Bs_0}{2}\right)
       +\frac32+\frac{\log K_0+2}{K_0},&&\\
 c_0&=(\sigma\ell-1)\log\rho-b_0.&&&
\end{align*}
The two numerical quantities to be certified are
\begin{align}
 F={}&\kappa c_0-\frac{2\log(\kappa_+\ell w_0)+3}{w_0}
       -g_0\ell\bigl(ra+\kappa_+s_0a_2\bigr),\label{direct-F}\\
 Q={}&\frac12\left(\mu\kappa_+\ell\log\rho+
          \frac{\log(10M_0w_0)+1}{w_0}\right),\label{direct-Q}\\
 M_0={}&\max\left\{\frac{\ell\kappa_+s_0}{8},
                    \frac{\ell r}{382w_0}\right\}.\notag
\end{align}

\begin{lemma}\label{direct-lemma}
For each row of Table~\ref{direct-table}, every quintic solution with
$191\le n\le B$ and $u>10^{1000}$ satisfies $n<Q$.
\end{lemma}

\begin{table}[htbp]
\centering
\caption{Direct Laurent certificates with $\mu=3/5$.
The bounds for $F$ are rounded down and those for $Q$ up.}
\label{direct-table}
\begin{tabular}{rrrrrrr}
\toprule
$B$ & $\kappa$ & $\ell$ & $r$ & $\rho$ & $F>$ & $Q<$\\
\midrule
271 & 42 & 6 & 22 & 16 & 0.668 & 209.613\\
199 & 44 & 5 & 23 & 19 & 1.215 & 194.338\\
193 & $433/10$ & 5 & 23 & 19 & 0.081 & 191.247\\
191 & $216/5$ & 5 & 23 & 19 & 0.176 & 190.805\\
\bottomrule
\end{tabular}
\end{table}

\begin{proof}[Proof of Lemma~\ref{direct-lemma}]
Use $\alpha_1=u^4/v$, $\alpha_2=2$, $b_1=n$, $b_2=4$ and
$D=1$, as in \eqref{quintic-Lambda}. In Laurent's Theorem~1 choose
\[
 K=\lceil\kappa w\rceil,\quad L=\ell,\quad R_1=1,\quad S_1=\ell,
 \quad R_2=r,\quad
 S_2=\left\lfloor\frac{(K-1)\ell}{r}\right\rfloor+1.
\]
Thus $R=r$, $S=\ell+\lfloor(K-1)\ell/r\rfloor$ and $N=K\ell$.
The first cardinality condition holds since the powers
$2^0,2^1,\ldots,2^{\ell-1}$ are distinct.
For the second condition, the integers $4i+nj$, with
$0\le i<r$, $0\le j<S_2$, are distinct. Equality implies
$n\mid i-i'$ because $n$ is odd, and $|i-i'|<r<n$.
Their number is $rS_2>(K-1)\ell$.

Take $a_1=aw$. Since $\log\alpha_1<4\log2/n$ and
$h(\alpha_1)=4w$, both height conditions in Laurent's theorem hold.
Also,
\[
 \kappa w\le K\le\kappa_+w,\qquad K\ge K_0,
 \qquad \frac SK\le s_0,
 \qquad g:=\frac14-\frac{K\ell}{12rS}\le g_0.
\]
To control the factorial product in that theorem we use, for $K\ge2$,
\begin{equation}\label{factorial-bound}
 -\frac{2}{K(K-1)}\log\prod_{j=1}^{K-1}j!
 \le-\log K+\frac32+\frac{\log K+2}{K}.
\end{equation}
Indeed, $\log(j!)\ge j\log j-j$, and comparison of
$\sum_{j=1}^{K-1}j\log j$ with $\int_1^{K-1}t\log t\,dt$ gives
\[
 \log\prod_{j=1}^{K-1}j!\ge
 \frac{(K-1)^2}{2}\log(K-1)
 -\frac{(K-1)^2-1}{4}-\frac{K(K-1)}2.
\]
Now use $\log(K/(K-1))\le2/K$ to obtain
\eqref{factorial-bound}. Consequently the quantity called $b$ in
Laurent's theorem satisfies $\log b\le b_0$.

The certified intervals give $b_0,c_0>0$ and $\log(\kappa_+\ell w_0)>1$.
The function $\log(\kappa_+\ell w)/w$ therefore decreases for
$w\ge w_0$. The left side of the theorem's numerical condition,
divided by $w$, is at least
\[
 \kappa c_0-\frac{2\log(\kappa_+\ell w_0)}{w_0}
      -g_0\ell(ra+\kappa_+s_0a_2).
\]
The right side is $\varepsilon(N)/w$. Since $N\ge10$, the estimate
$N!\le N^N$ gives
$\varepsilon(N)\le2+2\log(2N)/N<3$.
Thus the certified inequality $F>0$ verifies this condition for every
$w\ge w_0$.

Write
\[
 M=\max\{\ell S/8,\ell r/(2n)\}\le M_0w.
\]
The program also verifies
$1<\log(10M_0w_0)<382w_0$ and $382w_0>1$.
Together with \eqref{quintic-Lambda}, these show
$M\Lambda<10M_0w\exp(-382w)<1$ for $w\ge w_0$.
The modified form in Laurent's theorem has modulus at most
$M\Lambda e^{M\Lambda}$, and its modulus exceeds
$\rho^{-\mu K\ell}$. Since $M\Lambda<1$, taking logarithms gives
\[
 2nw<\mu K\ell\log\rho+\log(10M_0w)+1.
\]
Finally, $\log(10M_0w)/w$ decreases on $w\ge w_0$.
Using $K/w\le\kappa_+$ proves $n<Q$.
\end{proof}

Starting with $B=271$, the four rows of Table~\ref{direct-table}
give successive prime upper bounds $199$, $193$, $191$ and $181$.
Together with \eqref{n-coprime}, this yields \eqref{E}. Since
$2x+1=u^n>10^{1000n}\ge10^{23000}$, the lower bound in \eqref{huge}
follows as well. This proves Theorem~\ref{quintic-main}.

\FloatBarrier
\section{Counting the quintic solutions}\label{sec:counting}
\subsection{The gap between two solutions}

\begin{proof}[Proof of Proposition~\ref{gap}]
For $i=1,2$, inequality \eqref{CF-ineq} gives
\[
 0<\frac{v_i}{u_i^4}-2^{-4/n}<u_i^{-2n}.
\]
The two fractions are reduced and have distinct denominators. Their
difference in absolute value is therefore at least $1/(u_1^4u_2^4)$.
Since both lie on the same side of $2^{-4/n}$, this difference is
strictly less than $u_1^{-2n}$. Hence
$u_2^4>u_1^{2n-4}$, as required. The last assertion follows from
$u_1>10^{1000}$ and $n\ge23$.
\end{proof}

\subsection{An upper bound for the number of solutions}

For a reduced rational number $p/q$, with $q>0$, write
$H(p/q)=\max\{|p|,q\}$. If the primitive irreducible polynomial
of an algebraic number $\alpha$ is
$a\prod_{i=1}^d(X-\alpha_i)\in\Z[X]$, its Mahler measure is
\[
 M(\alpha)=|a|\prod_{i=1}^d\max\{1,|\alpha_i|\}.
\]
We use the following real specialization of the
principle of Bombieri and Mueller \cite[Section~II]{BM} in the form
given by Mosunov \cite[Lemma~2.9 of the preprint]{Mosunov}, taking
$\alpha_1=\alpha_2=\alpha$. The bounds on $\tau$ below already imply
the lower bound on $t$ in that lemma, since
$(t-2/d)^2>2-dt^2$ and $t>2/d$.

\begin{lemma}[Thue--Siegel principle]\label{thue-siegel}
Let $\alpha$ be a real algebraic number of degree $d\ge3$. Choose
real numbers $t,\tau$ such that
\[
 0<t<\sqrt{2/d},\qquad
 \sqrt{2-dt^2}<\tau<t-2/d,
\]
and put
\[
 \lambda=\frac2{t-\tau},\qquad
 A=\frac{t^2}{2-dt^2}\left(\log M(\alpha)+\frac d2\right),\qquad
 \delta=\frac{dt^2+\tau^2-2}{d-1}.
\]
If two reduced rational numbers $p_i/q_i$, $q_i>0$, satisfy
\[
 \left|\alpha-\frac{p_i}{q_i}\right|
 <\bigl(4e^A H(p_i/q_i)\bigr)^{-\lambda}
 \qquad(i=1,2),
\]
then
\begin{equation}\label{thue-siegel-bound}
 \log H(p_2/q_2)+A+\log4
 \le\delta^{-1}\bigl(\log H(p_1/q_1)+A+\log4\bigr).
\end{equation}
\end{lemma}

\begin{proof}[Proof of Theorem~\ref{quintic-count}]
Fix $n\in\mathcal E$ and set
\[
 \alpha=2^{-4/n},\qquad \mu=n/2,\qquad g=(n-2)/2.
\]
For a quintic solution put $p=v$ and $q=u^4$. Then $p/q$ is reduced,
and \eqref{CF-ineq} gives
\begin{equation}\label{count-approx}
 0<p/q-\alpha<q^{-\mu},\qquad q>10^{4000}.
\end{equation}
Since $v<u^4$, as shown in Section~\ref{sec:quintic}, we have $H(p/q)=q$.
The number $\alpha$ has degree $n$. Indeed, if $\beta=2^{1/n}$,
then $\alpha=\beta^{-4}$. Choose integers $a,b$ with $-4a+nb=1$.
The identity $\beta=\alpha^a2^b$ gives
$\mathbb Q(\alpha)=\mathbb Q(\beta)$, and $X^n-2$ is Eisenstein at $2$. Consequently the primitive irreducible
polynomial of $\alpha$ is $16X^n-1$, and $M(\alpha)=16$.

Order any finite collection of solutions by increasing $u$, or
equivalently by increasing $q$, and write $\ell_i=\log q_i$.
Proposition~\ref{gap} implies
\begin{equation}\label{count-log-gap}
 \ell_s>g^{s-r}\ell_r\quad(s>r),\qquad
 \ell_i>4000\log10\,g^{i-1}>9200g^{i-1}.
\end{equation}

Choose $a,b$ and $m$ from Table~\ref{count-table} and put
\begin{align*}
 t&=\sqrt{\frac2{n+a^2}},& \tau&=bt,&
 \lambda&=\frac{\sqrt{2(n+a^2)}}{1-b},\\
 A&=\frac{\log16+n/2}{a^2},&
 \delta&=\frac{2(b^2-a^2)}{(n+a^2)(n-1)},&
 L&=9200.
\end{align*}
For exact verification replace $A+\log4$ by the rational upper bound
\[
 \widehat B=\frac{n/2+14/5}{a^2}+\frac75,
\]
using $\log16<14/5$ and $\log4<7/5$. Define the rational quantities
\begin{equation}\label{count-checks}
 U=\frac{\lambda^2(1+\widehat B/L)^2}{\mu^2},\qquad
 V=\frac{1+\widehat B/L}{\delta g^m}.
\end{equation}
Every row satisfies $0<a<b<1$, $U<1$ and $V<1$.
The displayed upper bounds are rounded upwards.

\begin{table}[htbp]
\centering
\caption{Parameters giving $N_n\le m$. The last two columns give
strict rational upper bounds for $U$ and $V$.}
\label{count-table}
\begin{tabular}{rrrrrr}
\toprule
$n$ & $a$ & $b$ & $m$ & $U<$ & $V<$\\
\midrule
23 & $17/100$ & $37/100$ & 4 & 0.9747 & 0.2034\\
29 & $9/50$ & $11/25$ & 4 & 0.9862 & 0.0804\\
$31\le n\le181$ & $1/5$ & $9/20$ & 3 & 0.9416 & 0.9868\\
\bottomrule
\end{tabular}
\end{table}

For the first two rows the bounds follow by substituting the
displayed rational numbers into \eqref{count-checks}. For the last
row, $1+\widehat B/L=(92714+125n)/92000$.
Apart from positive constant factors, $U$ and $V$ are respectively
\[
 \frac{(n+1/25)(92714+125n)^2}{n^2},\qquad
 \frac{(n+1/25)(n-1)(92714+125n)}{(n-2)^3}.
\]
The signs of their derivatives on $[31,181]$ are the signs of
\[
 125n^2-92714n-\frac{185428}{25},\qquad
 -93344n^2-\frac{4808878}{25}n+\frac{4728664}{25},
\]
both negative. Hence it suffices to substitute $n=31$, which gives
the last two entries of the table. All these verifications use
rational arithmetic only.

We next check the hypotheses of Lemma~\ref{thue-siegel} with $d=n$.
We have $0<t<\sqrt{2/n}$ and
$\sqrt{2-nt^2}=at<bt=\tau$.
Since $U<1$, we also have $\lambda<\mu<n$, so
$\tau<t-2/n$. For every solution, \eqref{count-approx} gives
$\ell_i>4000\log10>L$ and therefore
\[
 \lambda\left(1+\frac{A+\log4}{\ell_i}\right)
 <\lambda(1+\widehat B/L)<\mu.
\]
Together with \eqref{count-approx}, this proves the approximation
hypothesis of that lemma for every such solution.

If there were $m+1$ solutions, Lemma~\ref{thue-siegel} would apply
to the first and last in order of increasing $q$. Since $A+\log4>0$,
it would give
\[
 \frac{\ell_{m+1}}{\ell_1}
 <\delta^{-1}\left(1+\frac{\widehat B}{L}\right)=Vg^m<g^m,
\]
whereas \eqref{count-log-gap} gives $\ell_{m+1}>g^m\ell_1$.
Thus $N_n\le m$. There are 32 primes between $31$ and $181$,
inclusive, so the total is at most $2\cdot4+32\cdot3=104$.
By Theorem~\ref{quintic-main} this counts all integer exponents $n\ge3$.
\end{proof}

\begin{proof}[Proof of Corollary~\ref{quintic-distinct-x}]
Fix $x$ and write $T=x^5+(x+1)^5$. By Theorem~\ref{quintic-main}, every admissible exponent is prime. If two distinct primes $p,q$
were exponents for $T$, then $pq\mid v_{\ell}(T)$ for every prime
$\ell\mid T$. Thus $T$ would be a $pq$th power, contradicting the same
theorem. For a given exponent the positive base is unique.
\end{proof}

\FloatBarrier
\section{Computational verification}\label{sec:computations}

The two verification programs, their certificates, and the execution reports
are available at \url{https://github.com/apinter1967/consecutive-powers}.
To reproduce the computations, extract the
\href{https://github.com/apinter1967/consecutive-powers/blob/eeaf345b2ad9474499ac401ad8ff262ffd0caf90/consecutive_powers_verification.zip}{ZIP archive at commit \texttt{eeaf345}}.
With Python 3.11 or later, run
\begin{quote}\small\ttfamily
python verify\_consecutive\_powers.py {-}{-}outdir checks\\
python verify\_solution\_count.py {-}{-}outdir checks
\end{quote}
The first program verifies the numerical bounds in
Tables~\ref{bounds} and~\ref{direct-table}, tests the 83 primes
in the congruence sieve, and computes the continued fractions for
the 58 primes up to $317$ in the stated range. Logarithms are
enclosed by range reduction and the series for
$2\operatorname{arctanh}z$, with a bound on its positive tail.
The program uses decimal interval arithmetic with directed rounding
and integer square roots to enclose radicals. The endpoints used for
continued fractions are certified by \eqref{alpha-cert}.
The second program uses rational arithmetic to check
Table~\ref{count-table}, the derivative bounds on $[31,181]$ and the
sum of the 34 individual bounds. The reports record the software
version and source checksums.

\section*{Acknowledgements}
The author used ChatGPT (OpenAI) to polish the mathematical English and
to assist with the computational parts of this work, including the
development and checking of the verification programs. The author
takes full responsibility for the content of the paper.

\end{document}